\documentclass[12pt]{amsart}
\newtheorem{thm}{Theorem}
\newtheorem{lem}{Lemma}
\newtheorem{cor}{Corollary}
\newtheorem{prop}{Proposition}

\newcommand{\nd}{|\!\!/}
\newcommand{\subs}{\substack}

\begin{document}
%%% \topmatter
\title[Counting squarefree discriminants]{Counting squarefree 
discriminants of trinomials under abc}
%% \endtitle
%\vskip 3mm
\author{Anirban Mukhopadhyay, M. Ram Murty and Kotyada Srinivas}
\address{Institute of Mathematical Sciences,
CIT Campus, Tharamani, Chennai 600 113, India}
\email[Anirban Mukhopadhyay]{anirban@imsc.res.in}
\email[Kotyada Srinivas]{srini@imsc.res.in}
\address{Department of Mathematics and Statistics,
Jeffery Hall, Queen's University,
Kingston, Ontario, K7L 3N6, Canada}
\email[M. Ram Murty]{murty@mast.queensu.ca}


\begin{abstract}
For an odd positive integer $n\ge 5$, assuming the truth of the $abc$
conjecture, we show that for a positive proportion of pairs 
$(a,b)$ of integers the 
trinomials of the form $t^n+at+b \ \ (a,b\in \mathbb Z)$
are irreducible and their discriminants are squarefree.
\end{abstract}

\maketitle

\section{Introduction}
Let $D_f$ be the discriminant of the trinomial 
\begin{equation}\label{trino}
f(t)=t^n+at+b \ \ (a,b\in \mathbb Z)
\end{equation} 
where $\mathbb Z$ denotes the set of integers.
For positive integers $A>1, B>1$ we define
$\mathcal M_n(A,B)$ to be the set of $(a,b)$ with
$A\le |a| \le 2A$, $B\le |b|\le 2B$ such that
$f(t)$ is irreducible and $D_f$ is squarefree.
It is reasonable to expect that for $A$, $B$ tending to infinity,
$$M_n(A,B) \sim c_n AB, $$
for some positive constant $c_n$.  This is probably very difficult
to prove.  We will apply the $abc$ conjecture to show
that $M_n(A,B) \gg AB$.  Recall that
the $abc$-conjecture, first formulated in 1985 by Oesterl\'e
and Masser is the following statement.

Fix $\epsilon>0$. If $a,b$ and $c$ are coprime positive integers
satisfying $a+b=c$, then
$$c\ll_{\epsilon} N(a,b,c)^{1+\epsilon}$$  
where $N(a,b,c)$ is the product of dictinct primes dividing $abc$.\\

Our main theorem is as follows
\begin{thm}\label{main}
Assume the truth of the $abc$ conjecture.
Let $n\ge 5$ be odd and $n\equiv 1\mod 4$. 
Let $A$ be sufficiently large and $B>A^{1+\delta_0}$ for some fixed $\delta_0>0$.
Then
$$\mathcal M_n(A,B)\gg AB$$
where the implied constants may depend on $n$. 
\end{thm}

\noindent {\bf Remark.}  The cases $n=2$ and $n=3$ of the theorem
can be treated without the use of the $abc$ conjecture.  Indeed,
the case $n=2$ reduces to counting the number of $a,b$ with $a^2 - 4b$
squarefree. This question is answered in \cite{murty-class}
as Theorem 3.  The case $n=3$ can be dealt with along the same lines.
Indeed, First, one counts the number of such pairs $(a,b)$
such that $4a^3+27b^2$ is squarefree.  This is easily done by
fixing $a$, using Theorem 3 of \cite{murty-class} and
then summing over $a$.  A cognate result is derived in \cite{chakra-murty}.
The case $n=4$ can be treated using the simple asymptotic sieve
as in \cite{hooley}.  In this case, we essentially need to
count how ofter $27a^4 + 256b^3 $ is squarefree.  Fixing $a$,
we are reduced to determining how often the value of a cubic polynomial
is squarefree.  Following the method of Chapter 4 of
\cite{hooley}, we easily derive the required result.  An appropriate
modification of this leads to an answer to the question under
consideration.  We leave the details to the reader.

Now we describe an application of the theorem.
In \cite{osada}, Osada showed that the Galois group of (\ref{trino})
is isomorphic to $S_n$ provided
\begin{enumerate}
\item $f(t)$ is irreducible over $\mathbb Q$,
\item $((n-1)a,nb)=1$.
\end{enumerate}
Moreover, if $K_f$ is the splitting field of $f(t)$ over $\mathbb Q$,
then $K_f$ is unramified at all finite primes over 
$\mathbb Q(\sqrt D_f)$ with the alternating
group $A_n$ of degree $n$ as the Galois group.

Using theorem \ref{main}, we prove the following quantitative version of
Osada's result.

\begin{cor}\label{unrami}
Assume the truth of $abc$ conjecture.
Let $n\ge 5$ be odd and $n\equiv 1\mod 4$. 
Also, let $\mathcal N_n(X)$ be the number 
of quadratic number fields of the form $\mathbb Q(\sqrt {D_f})$ 
with $|D_f|\le X$ which has a Galois extension with
Galois group $A_n$ and unramified at all finite primes. 
Then for large $X$,
$$\mathcal N_n(X)\gg X^{\frac{1}{n}+\frac{1}{n-1}}$$
where the implied constant may depend on $n$. 
\end{cor}

In order to prove the theorem we need to count irreducible polynomials
with square free discriminants. In section $2$, we show that almost all polynomials of
the specific form under consideration are irreducible. In section $3$, we show that
a positive proportion of the polynomials have square free discriminants. Sections
$4$ and $5$ provide the technical details needed in section $3$. The last section contains
the conclusion of the proof.

\section{Counting irreducible polynomials}

We start with a result due to S. D. Cohen \cite{sdc}
regarding the number of 
irreducible polynomials of a certain form over finite fields. Before
stating it we need to introduce some notations.
For a fixed prime $p$, let $g(t), h(t)$  be monic, relatively 
prime polynomials in $\mathbb F_p[t]$ satisfying
$$n=\deg g>\deg h\ge 0$$
and
$$g(t)/h(t)\neq g_1(t^p)/h_1(t^p), \ 
              {\rm for \ any} \ g_1(t), h_1(t)\in \mathbb F_p[t].$$
Let $L$ be the spliting field of 
$Q(y)=g(t)-yh(t)\in\mathbb F_p[t][y]$ over $\mathbb F_p(y)$
and $G$ be its Galois group. 
Let $\mathbb F_{p^f}$ be the 
maximal algebraic extension of $\mathbb F_p$ in $L$.
For any $\sigma\in G$, let $L_\sigma$ denote
the subfield of $L$ fixed by $\sigma$. We define 
$$G^*=\left\{\sigma\in G| L_\sigma\cap \mathbb F_{p^f} =\mathbb F_p \right\}.$$
We consider $G$ to be a subgroup of $S_n$. Let
$G_n=\{ \sigma\in G| \sigma \ {\rm is \ a} \ n-{\rm cycle}\}$ and
$G_n^*=G^*\cap G_n$. 
We define $\pi(g,h)$ to be the number of irreducible polynomials of
the form $hP(g/h)$, where $P$ is a linear monic polynomial in 
$\mathbb F_p[t]$.
Now we state a particular case of Theorem 3 in \cite{sdc}.

\begin{lem}\label{cohen}
$$\pi(g,h)=\frac{|G_n^*|}{|G^*|}p+O(\sqrt p).$$
\end{lem}

For a fixed $a\in\mathbb F_p$, let $g_a(t)=t^n+at$.
For $g=g_a$ and $h=1$, we get from \cite{bsd}, that
$$G=G^*=S_n$$
whenever $(p,2n(n-1))=1$.
Also 
$$\frac{|G_n^*|}{|G^*|}=\frac{1}{n}.$$
Hence from Lemma \ref{cohen}, we have
\begin{equation}\label{pi}
\pi(g_a,1)=\frac{p}{n}+O(\sqrt p).
\end{equation}
Clearly, for a fixed $a\in\mathbb F_p$, 
$\pi(g_a,1)$ is the number of irreducible polynomials of 
the form $t^n+at+b$ with $b\in\mathbb F_p$.

For a prime $p$ we define $\mathcal S_p$ and $\mathcal T_p$ as follows
$$\mathcal S_p=\left\{(a,b)\in (\mathbb F_p)^2| ~t^n+at+b~~ 
  {\rm is \ reducible \  over}~~\mathbb F_p \right\}$$

$$\mathcal T_p=\left\{(a,b)\in (\mathbb F_p)^2| ~t^n+at+b~~ 
  {\rm is \ irreducible \  over}~~\mathbb F_p \right\}$$
and let $s_p=|\mathcal S_p|$ and
$t_p=|\mathcal T_p|$. From (\ref{pi}),
varying over $a\in\mathbb F_p$, we get the following lemma
estimating $t_p$.

\begin{lem}\label{p-irr}
If $p$ does not divide $2n(n-1)$, then
$$t_p=\frac{p^2}{n}+O(p^{3/2}).$$
\end{lem}

Now we introduce the following notations.
$$\mathcal T(A,B)=\left\{(a,b)\in \mathbb Z^2| ~t^n+at+b~~ 
  {\rm is \ irreducible}, ~~A\le |a| \le 2A ,B\le |b|\le 2B\right\}.$$

The proof of the following proposition estimating the cardinality of 
$\mathcal T(A,B)$, 
closely follows the method outlined in exercise no.12,
page 169 of \cite{bourbaki}.
\begin{prop}\label{irr}
$$|\mathcal T(A,B)|=AB+
                o\left( AB \right).$$
\end{prop}

\proof 
We observe that $s_p+t_p=p^2$. So, from lemma \ref{p-irr}, we get
$$s_p=p^2\left(1-\frac{1}{n}\right)+O(p^{3/2}).$$
For a square free integer $d$, let
$$\phi_d:\mathbb Z^2\rightarrow (\mathbb Z/d\mathbb Z)^2$$
be the reduction modulo $d$. 
Let
$$H\subset\left\{(a,b)\in\mathbb Z^2|\ A \le |a| \le 2A,\ 
                                       B \le |b| \le 2B \right\}$$
and $H_d$ be the image of $H$ under $\phi_d$.
The number of elements of $H$ which are mapped to the same 
element of $(\mathbb Z/d\mathbb Z)^2$ under $\phi_d$ does not
exceed $([2A/d]+1)([2B/d]+1)$. We deduce
\begin{eqnarray*}
|H| & \le & |H_d|([2A/d]+1)([2B/d]+1)\\
    & \ll & \frac{|H_d|}{d^2}AB \\
    & \ll & \left(\prod_{p|d}\frac{|H_p|}{p^2}\right)AB.
\end{eqnarray*}
Now we set
$$H=\left\{(a,b)\in \mathbb Z^2| ~t^n+at+b~~ 
  {\rm is \ reducible, with}~~A\le |a|\le 2A,B\le |b|\le 2B \right\}.$$
Then $H_p\subset \mathcal S_p$ for each prime $p$.
From above, we have
$$|H|\ll \left(\prod_{p|d}\frac{s_p}{p^2}\right)AB
     \ll \prod_{p|d}\left(1-\frac{1}{n}\right)AB.$$

For $\epsilon>0$, we choose $m>1$ such that
$$\left(1-\frac{1}{n}\right)^m<\epsilon.$$
Let $p_1,p_2,\cdots, p_m$ be the first $m$ primes not dividing $2n(n-1)$.
By choosing $d=p_1\cdots p_m$, we get
$$|H|\ll \epsilon AB.$$
Hence the proposition follows.
$\hfill{\Box}$

Corollary \ref{unrami} is the quantitative version of the following 
result due to Osada (see corollary 2, \cite{osada}).

Let $K_f$ be the spliting field of $f(t)$ over $\mathbb Q$
and $G$ be the Galois group $Gal(K_f/\mathbb Q)$.
\begin{lem}\label{osada}
Let $f(t)=t^n+at+b$ be a polynomial in $\mathbb Z [t]$,
where $a=a_0c^n$ and $b=b_0c^n$ for some integer $c$.
Then the Galois group $G$ is isomorphic to $S_n$
if the folowing conditions are satisfied.
\begin{itemize}
\item[(1)] $f(t)$ is irreducible over $\mathbb Q$.
\item[(2)] $(a_0c(n-1),nb_0)=1$.
\end{itemize}
Moreover, $K/\mathbb Q(\sqrt{D_f})$ is unramified at all finite places. 
\end{lem}

\section{Counting square-free discriminants}

Let $T(a,b)=(n-1)^{n-1}a^n+n^nb^{n-1}$ for integers $a,b$.
For $n\equiv 1 \mod 4$, we observe that discriminant  $D_f=T(a,b)$.
For sufficiently large positive real numbers $A,B$,
Let $D(A,B)$ be the number of square free integers $d$ 
with at least one solution to 
\begin{equation}\label{eqd}
d=T(a,b)\ {\rm where}\ 
             A\le |a|\le 2A, B\le |b|\le 2B
\end{equation}
and $((n-1)a,nb)=1$. Using ideas from \cite{sound}
we now find a lower bound for $D(A,B)$.  
For a square free number $d$, let $R_0(d)$ denote the number of
solutions to (\ref{eqd}). We have 

\begin{lem}\label{r0d}
$$\sum_{d}R_0(d)\gg AB.$$
\end{lem}
\begin{lem}\label{r0dsq}
$$\sum_{d}R_0(d)^2\ll AB.$$
\end{lem}

Proof of these two lemmas will be presented in the next section.
Assuming them, we are ready to prove the following result
giving a lower bound for $D(A,B)$.

\begin{prop}\label{sqfree}
$$D(A,B)\gg AB.$$
\end{prop}

\proof
By Cauchy-Schwarz inequality,
$$D(A,B)\ge \left(\sum_{d}R_0(d)\right)^2
                     \left(\sum_{d}R_0(d)^2\right)^{-1}.$$
Hence the result follows by lemmas \ref{r0d} and \ref{r0dsq}.
$\hfill{\Box}$

\section{Proof of Lemma \ref{r0d}} 

We define a new polynomial $H(a,b)=T(a,b)T(-a,b)$.
Let $\mathcal M_1$ be the set of pairs $(a,b)$ of integers with 
$A\le a\le 2A$ and $B\le |b|\le 2B$ such that 
$H(a,b)$ is not divisible by
the square of any prime $p\le \log B$. We put 
$M_1=\#\mathcal M_1$ and $P=\prod_{p\le \log B}p$. 
We observe that
$\sum_{l^2|(\alpha, P^2)}\mu(l)=1$,
or $0$ depending on whether $p^2\nd\alpha$ for all $p\le\log B$, or not.
Thus 
\begin{equation}\label{sumb}
M_1=\sum_{A\le a\le 2A}\sum_{B\le |b|\le 2B}\sum_{l^2|(H(a,b), P^2)}\mu(l)=
    \sum_{A\le a\le 2A} \sum_{l|P}\mu(l)\sum_{\subs{B\le |b|\le 2B \\ 
H(a,b)\equiv 0\mod l^2}}1.
\end{equation}
Let 
$$\rho_a(p)=|\{b\mod p|H(a,b)\equiv 0\mod p\}|.$$

Clearly $\rho_a(l)$ is a multiplicative function of $l$.
For a prime $p\nd an(n-1)$ and an integer $\alpha\ge 1$,
$$\rho_a (p^{\alpha})=\rho_a (p)\le n.$$
We divide the sum over $b$ in (\ref{sumb}) into intervals of length
$l^2$. We see that this sum is 
$$2B\rho_a(l^2)/l^2+O(\rho_a(l^2))$$
Thus,

\begin{eqnarray*}
&&\sum_{l|P}\mu(l) \sum_{\subs{B\le |b|\le 2B \\ 
H(a,b)\equiv 0\mod l^2}}1 \\
&=&2B \sum_{l|P}\mu(l)\frac{\rho_a(l^2)}{l^2}+
     O(\sum_{l|P}\mu(l)\rho_a(l^2))\\
    &=& 2B\prod_{p|P}\left(1-\frac{\rho_a(p)}{p^2}\right)
                                 +O(\tau(P))\\
    &=& 2B\exp\left(-\sum_p\frac{\rho_a(p)}{p^2}
         +O\left(\sum_{p>\log B}\frac{1}{p^2}\right)\right)
        +O(B^{\epsilon}), \\
    &=& 2cB+o(B)\ (\rm{for \ some \ constant}\ c>0) 
\end{eqnarray*}
where $\tau(\alpha)$ denote the divisor function and
 we use the observation that $P\asymp B$.
Summing over all choices of $a$ we get from (\ref{sumb})
$$M_1= 2cAB+o(AB)\ (\rm{for \ some \ constant}\ c>0).$$ 

Let $\mathcal M_2$ be the set of pairs $(a,b)$, $A\le a\le 2A$ and 
$B\le |b|\le 2B$ such that $H(a,b)$ is divisible by square of a prime
$p\in (\log B, B]$. Also let
$M_2=\#\mathcal M_2$. Then
\begin{eqnarray*}
M_2 &=&\sum_{A\le a\le 2A}\sum_{\log B<p\le B} 
                 \sum_{\subs{B\le |b|\le 2B\\H(a,b)\equiv 0\mod p^2}}1\\
    &=&2B\sum_{A\le a\le 2A}\sum_{\log B<p\le B}\frac{\rho_a(p^2)}{p^2}
       +O\left(\sum_{A\le a\le 2A}\sum_{\log B<p\le B}\rho_a(p^2)\right)\\
\end{eqnarray*}
The first term is
$$\ll AB\sum_{p>\log B}\frac{1}{p^2}\ll \frac{AB}{\log B}=o(AB).$$
The $O$-term is estimated as
$$\ll A \sum_{\log B<p\le B}1
  \ll \frac{AB}{\log B}=o(AB).$$
Then $\mathcal M_1 \setminus \mathcal M_2$ is the set of pairs $(a,b)$, 
$A\le a\le 2A$ and
$B\le |b|\le 2B$ such that both $T(a,b)$ and $T(-a,b)$ are not
divisible by square of a prime $p\le B$. 
We observe that $\#(\mathcal M_1\setminus \mathcal M_2)\ge M_1-M_2$.

We call a pair $(a,b)$ ``good'' if $T(a,b)$ is not
divisible by square of a prime $p > B$,
otherwise $(a,b)$ is called ``bad''. 

Now we claim that $(a,b)$ and $(-a,b)$ both cannot be bad.
Suppose both are bad. Then there are primes $p,q > B$
such that
$$T(a,b)=p^2d_1, \ \ T(-a,b)=q^2d_2.$$ 
Since $n$ is odd we get by multiplying
\begin{equation}\label{abc}
T(a,b)T(-a,b)=n^{2n}b^{2(n-1)}-(n-1)^{2(n-1)}a^{2n}=p^2q^2d_1d_2.
\end{equation}
If $p=q$, then $p$ divides $T(a,b)+T(-a,b)$ implying $p\le B$,
a contradiction.
Thus $p$, $q$ are distinct. 
Using $abc$-conjecture on the equation (\ref{abc}), we get
for any $\epsilon>0$
$$pq\ll (AB)^{1+\epsilon}$$
which is a contradiction as $p,q$ both are $>B$ and $B>A^{1+\delta_0}$
for a fixed $\delta_0>0$.

Hence among the pairs in $\mathcal M_1 \setminus \mathcal M_2$ half of them are good, 
and hence square-free as they are not divisible by square of a prime $\le B$.
Thus
$$\sum_{d\le X}R_0(d)\ge \frac{1}{2}(M_1-M_2)\gg AB.$$ 
This completes the proof.\\

\section{Proof of Lemma \ref{r0dsq}}

Let $a_1,a_2$ are in $[-2A,-A]\cup[A,2A]$ and 
$b_1,b_2$ are in $[-2B,-B]\cup[B,2B]$. Then $\sum_{d\le X}R_0(d)(R_0(d)-1)$
is bounded by the number of $(a_1,a_2,b_1,b_2)$ with 
$(a_1,a_2)\neq (b_1,b_2)$ and $T(a_1,b_1)=T(a_2,b_2)$.
Then 
$$(n-1)^{n-1}a_1^n+n^nb_1^{n-1}=(n-1)^{n-1}a_2^n+n^nb_2^{n-1}$$
which implies
$$(n-1)^{n-1}(a_1^n-a_2^n)=n^n(b_2^{n-1}-b_1^{n-1}).$$
Thus, for fixed $(a_1,a_2)$, number of possible $b_1$ and $b_2$ 
is $\ll X^{\epsilon}$. Hence
$$\sum_{d\le X}R_0(d)(R_0(d)-1)\ll X^{\epsilon}A^2.$$
Therefore, we have
$$\sum_{d\le X}R_0(d)^2\ll X^{\epsilon}A^2+AB\ll AB$$
completing the proof of Lemma \ref{r0dsq}.

\section{Proof of the Theorem and the Corollary}
Let $\mathcal D(A,B)$ be the set of $(a,b)$'s chosen exactly one for 
each $d$ counted in $D(A,B)$. Thus $D(A,B)=|\mathcal D(A,B)|$. 
Clearly 
$$\mathcal T(A,B)\cap \mathcal D(A,B) \subset \mathcal M(A,B).$$
Hence the theorem follows from propositions $1$ and $2$.

The corollary is a direct consequence of the theorem with
$A=X^{1/n}/4n$ and $B=X^{1/(n-1)}/4n^2$ and Lemma \ref{osada}
with $c=1$.

\vspace{0.5cm}
\noindent {\bf Acknowledgements.} We thank D. Suryaramana for indicating the
proof of proposition 1. We also thank R. Balasubramanian for
some useful discussions.


\begin{thebibliography}{DPR}
%\bibitem{A-C} N. Ankeny and S. Chowla: On the divisibility of the
%class numbers of quadratic fields, {\sl Pacific Journal of Math.},
%no. {\bf5}(1995), 321--324.

\bibitem{bsd} B. J. Birch and H. P. F. Swinnerton-Dyer: 
Note on a problem of Chowla, 
{\sl Acta Arith}, Vol.{\bf 5} (1959), 417-423.

\bibitem{bourbaki} Bourbaki:  
Algebra, Chapter 5
{\sl Springer Verlag}.


%\bibitem{CL} H. Cohen and H. W. Lenstra Jr.: Heuristics on class
%groups of number fields, {\sl Springer Lecture Notes}, {\bf 1068} in
%Number Theory Noordwijkerhout 1983 Proceedings.

\bibitem{chakra-murty}  K. Chakraborty and M. Ram Murty, 
On the number of real quadratic fields with class number
divisible by $3$, {\sl Proceedings of the American Math. Society,
\bf 131} (2002), No. 1, 41-44.



\bibitem{sdc} S. D.  Cohen: The distribution of polynomials over
finite fields,  
{\sl Acta Arith.}, Vol.{\bf 17} (1970), 255-271.

%\bibitem{D-H} H. Davenport and H. Heilbronn: On the density of
%discriminants of cubic fields II, {\sl Proc. Royal Soc.}, A {\bf 322} 
%(1971), 405--420.

%\bibitem{dickson} Dickson: Theory of Numbers. 
%{\sl Chelsea Publishing Company}, (1952).

%\bibitem{granville} Andrew Granville : $ABC$ allows us to count squarefrees,
%{\sl IMRN}, Vol.{\bf 19}, (1998), 991--1009.

%\bibitem{friesen-wam}Christian Friesen and Paul van Wamelen: 
%Class number of real quadratic function fields,
%{\sl Acta. Arith.} {\bf 81} (1997), no. 1, 45-55.

%\bibitem{Ho} T. Honda: A few remarks on class numbers of imaginary
%quadratic fields, {\sl Osaka J. Math.}, {\bf 12} (1975), 19--21.

\bibitem{hooley}  C. Hooley, Applications of sieve methods, 
Cambridge University Press, 1976.

\bibitem{murty-class} M. Ram Murty: Exponents of class groups of quadratic
fields, {\sl Topics in Number Theory (University Park, PA, 1997),
Math. Appl. {\bf 467}, Kluwer
Acad. Publ., Dordrecht}, (1999), 229--239.

\bibitem{osada} H. Osada: The Galois groups of polynomials 
$x^n+ax^l+b$, 
{\sl J. Number Theory}, {\bf 25} (1987), 230--238.



%\bibitem{rosen} Michael Rosen: Number Theory in Function Fields,
%{\sl Graduate Texts in Mathematics}, {\bf 210}. Springer-Verlag.
 
\bibitem{sound} K. Soundararajan: Divisibility of class numbers of
imaginary quadratic fields, {\sl J. London Math. Soc.} {\bf 61}
(2000), no. 2, 681--690.

%\bibitem{Na} T. Nagell: {\" U}ber die Klassenzahl imagin{\"a}r
%quadratischer
%Zahlkorper: {\sl Abh. Math. Sem. Univ. Hamburg}, {\bf 1} (1922),
%140--150.

%\bibitem{uchida} K. Uchida: Class numbers of cubic cyclic fields,
%{\sl J. Math. Soc. Japan}, Vol. {\bf 26}, no. 3, (1974), 447-453.


%\bibitem{Wi} P. Weinberger: Real quadratic fields with class numbers
%divisible by n, {\sl J. Number Theory}, {\bf 5} (1973), 237--241.

%\bibitem{Ya1} Y. Yamamoto: On unramified Galois extensions of
%quadratic number fields, {\sl Osaka J. Math.}, {\bf 7} (1970),
%57--76.
\end{thebibliography}
\end{document}